\documentclass[a4paper,12pt]{article}
\usepackage{amsmath, amssymb, amsthm}

\theoremstyle{plain}

\newtheorem{thm}{Theorem}[section]
\newtheorem{prop}[thm]{Proposition}
\newtheorem{lemma}[thm]{Lemma}

\newtheorem{corollary}[thm]{Corollary}

\newtheorem{question}[thm]{Question}

\theoremstyle{definition}

\begin{document}

\def\proof{\paragraph{Proof.}}
\def\endproof{\hfill$\square$}
\def\noproof{\hfill$\square$}

\def\Z{{\mathbb Z}}\def\N{{\mathbb N}} \def\C{{\mathbb C}}
\def\Q{{\mathbb Q}}\def\R{{\mathbb R}} \def\E{{\mathbb E}}
\def\P{{\mathbb P}}

\def\SS{{\cal S}}\def\TT{{\cal T}}\def\MM{{\cal M}}
\def\GG{{\cal G}}\def\Id{{\rm Id}}\def\Sym{{\rm Sym}}
\def\Supp{{\rm Supp}}\def\AA{{\cal A}}\def\UU{{\cal U}}
\def\BB{{\cal B}}\def\ord{{\rm ord}}

\title{\bf{On singular Artin monoids}}
 
\author{
\textsc{Eddy Godelle and Luis Paris}}

\date{\today}

\maketitle

\begin{abstract} 
In this paper we study some combinatorial aspects of the singular Artin monoids. Firstly, we show 
that a singular Artin monoid $SA$ can be presented as a semidirect product of a graph monoid with 
its associated Artin group $A$. Such a decomposition implies that a singular Artin monoid embeds 
in a group. Secondly, we give a solution to the word problem for the FC type singular Artin 
monoids. Afterwards, we show that FC type singular Artin monoids have the FRZ property. Briefly 
speaking, this property says that the centralizer in $SA$ of any non-zero power of a standard 
singular generator $\tau_s$ coincides with the centralizer of any non-zero power of the 
corresponding non-singular generator $\sigma_s$. Finally, we prove Birman's conjecture, namely, 
that the desingularization map $\eta: SA \to \Z [A]$ is injective, for right-angled singular 
Artin monoids. 
\end{abstract}


\section{Introduction}

Let $S$ be a finite set. Recall that a {\it Coxeter matrix} over $S$ is a matrix 
$M=(m_{s\,t})_{s,t\in S}$ indexed by the elements of $S$ and such that $m_{s\,s} =1$ for all $s 
\in S$, and $m_{s\,t} = m_{t\,s} \in \{2,3,4, \dots, +\infty\}$ for all $s,t \in S$, $s\neq t$. A 
Coxeter matrix $M=(m_{s\,t})$ is usually represented by its {\it Coxeter graph}, $\Gamma$, which 
is defined as follows. $S$ is the set of vertices of $\Gamma$, two vertices $s,t$ are joined by 
an edge if $m_{s\,t}\ge 3$, and this edge is labelled by $m_{s\,t}$ if $m_{s\,t}\ge 4$. If $a,b$ 
are two letters and $m \in \Z_{\ge 2}$, then we denote by $w(m:a,b)$ the word $\dots bab$ of 
length $m$. We take an abstract set $\SS=\{\sigma_s; s \in S\}$ in one-to-one correspondence with 
$S$, and we define the {\it Artin group} associated to $\Gamma$ to be the group $A=A_\Gamma$ 
presented by
$$
A=A_\Gamma= \langle \SS\ |\ w(m_{s\,t}: \sigma_s, \sigma_t) = w(m_{s\,t}: \sigma_t, \sigma_s) 
\text{ for } s,t \in S,\, s \neq t \text{ and } m_{s\,t} < +\infty \rangle\,.
$$
The {\it Coxeter group} associated to $\Gamma$ is the quotient $W=W_\Gamma$ of $A$ by the 
relations $\sigma_s^2=1$, $s \in S$.

Take $X \subset S$ and put $\SS_X= \{ \sigma_s; s \in X\}$. We denote by $A_X$ the subgroup of 
$A$ generated by $\SS_X$, and by $W_X$ the subgroup of $W$ generated by $\SS_X$. Let $\Gamma_X$ 
be the full subgraph of $\Gamma$ generated by $X$. Then $A_X$ is the Artin group associated to 
$\Gamma_X$ (see \cite{Lek} and \cite{Par}), and $W_X$ is the Coxeter group associated to 
$\Gamma_X$ (see \cite{Bou}). The subgroup $A_X$ is called {\it standard parabolic subgroup} of 
$A$, and $W_X$ is called {\it standard parabolic subgroup} of $W$.

We say that $\Gamma$ (or $A$) is of {\it spherical type} if $W$ is finite, that $\Gamma$ (or $A$) 
is {\it right-angled} if $m_{s\,t} \in \{2, +\infty\}$ for all $s,t \in S$, $s \neq t$, and that 
$\Gamma$ (or $A$) is of {\it type FC} if, for all $X \subset S$, either $W_X$ is finite, or there 
exist $s,t \in X$, $s \neq t$, such that $m_{s\,t} = +\infty$. Note that the spherical type Artin 
groups as well as the right-angled Artin groups are both FC type Artin groups. The number 
$n=|S|$ is called the {\it rank} of $A$.

The first (non-abelian) example of Artin group which has appeared in the literature is certainly 
the braid group 
$\BB_n$ introduced by Artin \cite{Art} in 1925. One of the most important works in the subject is 
a paper by Garside \cite{Gar} where the word problem and the conjugacy problem for $\BB_n$ are 
solved. Garside's ideas have been extended to all spherical type Artin groups by Brieskorn, Saito 
\cite{BrSa}, and Deligne \cite{Del} in 1972. These two papers, \cite{BrSa} and \cite{Del}, are 
the foundation of the theory of Artin groups, and, more specifically, of the spherical type Artin 
groups. Right-angled Artin groups are also known as {\it graph groups} or as {\it free partially 
commutative groups}. They have been widely studied, and their applications extend to various 
domains like parallel computation, random walks, and cohomology of groups. We mention, for 
example, the paper \cite{BeBr} where a group which is $FP_2$ but note finitely presented is 
constructed as a subgroup of some right-angled Artin group. Artin groups of type FC have been 
introduced by Charney and Davis \cite{ChDa} in 1995 in their study of the $K(\pi,1)$-problem for 
complements of infinite hyperplane arrangements associated to reflection groups.

In the same way as the braid group $\BB_n$ has been extended to the singular braid monoid 
$S\BB_n$ (see \cite{Bae} and \cite{Bir}), the Artin groups can be extended to the singular Artin 
monoids as follows. Take a new abstract set $\TT= \{ \tau_s; s \in S \}$ in one-to-one 
correspondence with $S$, and define the {\it singular Artin monoid} associated to $\Gamma$ to be 
the monoid $SA=SA_\Gamma$ presented as a monoid by the generating set $\SS \cup \SS^{-1} \cup 
\TT$, where $\SS^{-1} = \{\sigma_s^{-1}; s \in S\}$, and by the relations
$$
\begin{array}{cl}
\sigma_s \sigma_s^{-1} = \sigma_s^{-1} \sigma_s = 1\,, \quad &\text{for } s \in S\,,\\
\sigma_s \tau_s = \tau_s \sigma_s\,, \quad &\text{for } s \in S\,,\\
w(m_{s\,t}: \sigma_s, \sigma_t)= w(m_{s\,t}: \sigma_t, \sigma_s)\,, \quad &\text{for } s,t \in S, 
\, s \neq t, \text{ and } m_{s\,t} < +\infty\,,\\
w( m_{s\,t}-1: \sigma_s, \sigma_t) \tau_s = \tau_{s \wedge t} w(m_{s\,t}-1: \sigma_s, \sigma_t)\,, 
\quad &\text{for } s,t \in S, \, s \neq t, \text{ and } m_{s\,t} < +\infty\,,\\
\tau_s \tau_t = \tau_t \tau_s\,, \quad &\text{for } s,t \in S, \, s \neq t, \text{ and } m_{s\,t} 
=2\,,\\
\end{array}
$$ 
where $s \wedge t = s$ if $m_{s\,t}$ is even, and $s \wedge t = t$ if $m_{s\,t}$ is odd. Observe 
that we have an epimorphism $\theta: SA \to A$ which sends $\sigma_s^{\pm 1}$ to $\sigma_s^{\pm 
1}$ and $\tau_s$ to $\sigma_s$ for all $s \in S$, and that this epimorphism has a section $\iota: 
A \to SA$ which sends $\sigma_s^{\pm 1}$ to $\sigma_s^{\pm 1}$ for all $s \in S$. In particular, 
$A$ embeds in $SA$.

The combinatorial study of the singular Artin monoids (of spherical type) has been initiated by 
Corran \cite{Cor} with techniques inspired from \cite{BrSa}.

The purpose of the present paper is to study different combinatorial aspects of the monoid $SA$. 

Firstly, in Section 2, we prove that $SA$ can be decomposed as a semidirect product of a so-called 
graph monoid with the Artin group $A$. A consequence of this decomposition shall be that 
$SA$ embeds in a group.
Note that this last result has been previously proved
by Basset \cite{Bas} and Keyman \cite{Key} with completely different proofs. 

Sections 3 and 4 concern only singular Artin monoids of type FC. We solve the word problem for $SA$ in 
Section 3.
In \cite{FRZ}, Fenn, Rolfsen, and Zhu proved that the centralizer in the singular braid monoid of 
a standard singular generator $\tau_j$ is equal to the centralizer of any non-zero power of 
$\tau_j$, and that this centralizer coincides with the centralizer of any non-zero power of 
$\sigma_j$. This property, which we like to call {\it FRZ property}, is of importance in the 
study of singular braids, and, in particular, in the proof of Birman's conjecture for braid groups 
(see \cite{Par3}). In Section 4, we extend the FRZ property to all singular Artin monoids of type FC.

Define the {\it desingularization map} as the multiplicative homomorphism $\eta: SA \to \Z [A]$ 
which sends $\sigma_s^{\pm 1}$ to $\sigma_s^{\pm 1}$ and $\tau_s$ to $\sigma_s - \sigma_s^{-1}$ 
for all $s \in S$. One of the main questions in the subject, known as {\it Birman's conjecture}, 
is to determine whether the desingularization map is injective. This is known to be true for 
braid groups \cite{Par3} and for rank 2 Artin groups \cite{Eas}. We prove Birman's conjecture for 
right-angled Artin groups in Section 5.

The last section is dedicated to some questions for which we do not have any significant result 
but that deserve to be mentioned. 

\section{Semidirect product structure}

Our purpose in this section is to determine a decomposition of $SA$ as a semidirect product of a 
so-called graph monoid with the Artin group $A$. A first consequence of this decomposition 
is that a singular Artin monoid embeds in a group. This decomposition shall be also used in 
Sections 3 and 4 to solve the word problem and to prove the ``FRZ property'' for singular Artin monoids of 
type FC. 

Let $G$ be a (standard) graph, let $V$ be its set of vertices, and let $E=E(G)$ be its set of 
edges. Define the {\it graph monoid} of $G$ to be the monoid $\MM(G)$ 
presented as a monoid by
$$
\MM(G)= \langle V\ |\ uv=vu \text{ if } \{u,v\} \in E \rangle^+\,.
$$

The first (standard) graph that we shall consider is the graph $\hat \Omega$ defined by the following data.

\smallskip
$\bullet$ $\hat \Upsilon= \{ \alpha \tau_s \alpha^{-1}; \alpha \in A \text{ and } s \in S\}$ is 
the set of vertices of $\hat \Omega$;

\smallskip
$\bullet$ $\{ \hat u, \hat v\}$ is an edge of $\hat \Omega$ if $\hat u \hat v = \hat v \hat u$ in 
$SA$.

\begin{prop}
We have $SA= \MM( \hat \Omega) \rtimes A$.
\end{prop}

\proof
We have a homomorphism $f: \MM( \hat \Omega) \rtimes A \to SA$ defined by $f(\alpha)= \alpha \in 
SA$ for all $\alpha \in A$, and $f (\hat u)= \hat u \in SA$ for all $\hat u \in \hat \Upsilon$. 
Conversely, One can easily verify using the presentation of $SA$ that there is a homomorphism $g: 
SA \to \MM(\hat \Omega) \rtimes A$ which sends $\sigma_s^{\pm 1}$ to $\sigma_s^{\pm 1}\in A$, and 
$\tau_s$ to $\tau_s \in \hat \Upsilon$, for all $s \in S$. Obviously, $f \circ g = \Id$ and $g 
\circ f = \Id$.
\endproof

\begin{corollary}
$SA$ embeds in a group.
\end{corollary}

\proof
Let $\GG(\hat \Omega)$ be the group presented by
$
\GG(\hat \Omega)= \langle \hat \Upsilon\ |\ \hat u \hat v = \hat v \hat u \text{ if } \{\hat u, 
\hat v\} \in E(\hat \Omega) \rangle\,.
$
Then $\MM(\hat \Omega)$ embeds in $\GG(\hat \Omega)$ (see \cite{DuKr} and \cite{Dub}), thus $SA= 
\MM(\hat \Omega) \rtimes A$ embeds in $\GG (\hat \Omega) \rtimes A$.
\endproof 

\section{The word problem}

Note that a solution to the word problem for $SA$ will also give a solution to the word problem 
for the Artin group $A$. 
So, a reasonable approach would be to study the word problem for those singular Artin monoids whose
associated Artin groups have known solutions to the word problem.
In the case of Artin groups of type 
FC, a solution has been found by Altobelli \cite{Alt}. 
Another observation is that a given element of a graph monoid $\MM(G)$ has finitely many representatives
and these representatives can be easily listed.
Other solutions to the word 
problem for $\MM(G)$ can be found in \cite{CaFo}, \cite{Vie}, and \cite{DuKr}. Now, assume that $A$ is 
of type FC and consider the decomposition $SA= \MM (\hat \Omega) \rtimes A$ of the previous 
section. By the above observations, in order to solve the word problem for $SA$, it suffices to 
find an algorithm which decides whether two elements $\alpha \tau_s \alpha^{-1}$ and $\beta 
\tau_t \beta^{-1}$ of $\hat \Upsilon$ are equal, and, if not, whether they commute or not. Such 
an algorithm can be easily derived from Proposition 3.1 below together with Altobelli's solution 
to the word problem for $A$.

Define the graph $\Omega$ as follows.

\smallskip
$\bullet$ $\Upsilon= \{ \alpha \sigma_s \alpha^{-1}; \alpha \in A \text{ and } s \in S \}$ is the 
set of vertices of $\Omega$;

\smallskip
$\bullet$ $\{u,v\}$ is an edge of $\Omega$ if $uv=vu$ in $A$.

\begin{prop}
Assume $\Gamma$ to be of type FC. 
Then there exists an isomorphism $\varphi: \hat \Omega \to \Omega$ which sends $\alpha \tau_s 
\alpha^{-1}$ to $\alpha \sigma_s \alpha^{-1}$ for all $\alpha \in A$ and all $s \in S$.
\end{prop}

\bigskip
The remainder of the section is dedicated to the proof of Proposition 3.1.

Define the {\it Artin monoid} associated to $\Gamma$ to be the monoid $A^+= A_\Gamma^+$ 
presented as a monoid by
$$
A^+= \langle \SS\ |\ w(m_{s\,t}: \sigma_s, \sigma_t)= w(m_{s\,t}: \sigma_t, \sigma_s) \text{ for } 
s,t \in S,\, s \neq t, \text{ and } m_{s\,t}< +\infty \rangle^+\,.
$$
By \cite{Par2}, the natural homomorphism $A^+ \to A$ which sends $\sigma_s$ to $\sigma_s$ for all 
$s \in S$ is injective. 
We can define the {\it length function} $\lg: A^+ \to \N$ which associates to each element of
$A^+$ the length of any of its representatives with respect to the generating set $\SS$.
Since the defining relations of $A^+$ are homogeneous, this function is well-defined and is a 
homomorphism of monoids. For the same reason, we can define a partial order 
$\le_R$ on $A^+$ by setting $a \le_R b$ if there exists $c \in A^+$ such that $ca=b$.
Now, the 
following proposition is a mixture of several well-known facts on spherical type Artin groups.

\begin{prop}
Assume $\Gamma$ to be of spherical type.

\smallskip
(1) {\bf (Brieskorn-Saito \cite{BrSa}, Deligne \cite{Del})}. $(A^+, \le_R)$ is a lattice. The 
lattice operations of $(A^+,\le_R)$ are denoted by $\wedge_R$ and $\vee_R$.

\smallskip
(2) {\bf (Brieskorn-Saito \cite{BrSa}, Deligne \cite{Del})}. Let $s,t \in S$, $s \neq t$. Then 
$\sigma_s \vee_R \sigma_t = w(m_{s\,t}: \sigma_s, \sigma_t) = w(m_{s\,t}: \sigma_t, \sigma_s)$.

\smallskip
(3) {\bf (Brieskorn-Saito \cite{BrSa}, Deligne \cite{Del})}. Let $\Delta= \vee_R \{ \sigma_s; 
s\in S\}$. Then there exists a permutation $\mu: S \to S$ such that $\mu^2= \Id$ and $\Delta 
\sigma_s \Delta^{-1} = \sigma_{\mu(s)}$ for all $s \in S$.

\smallskip
(4) {\bf (Brieskorn-Saito \cite{BrSa}, Deligne \cite{Del})}. Each $\alpha \in A$ can be written as 
$\alpha= a \Delta^k$ with $a \in A^+$ and $k \in \Z$.

\smallskip
(5) {\bf (Charney \cite{Cha})}. Each $\alpha \in A$ can be uniquely written as $\alpha= ab^{-1}$ 
with $a,b \in A^+$ and $a \wedge_R b=1$. Such an expression $\alpha=ab^{-1}$ is called the {\rm 
Charney form} of $\alpha$.

\smallskip
(6) {\bf (Charney \cite{Cha})}. Let $\alpha \in A$ and $u,v \in A^+$ such that $\alpha=uv^{-1}$. 
Then $u=ac$ and $v=bc$, where $c=u\wedge_R v$ and $ab^{-1}$ is the Charney form of 
$\alpha$.
\noproof
\end{prop}

Let $s,t \in S$ and $\omega \in A$. We say that $\omega$ is an {\it elementary positive $(t,s)$-ribbon} if 
either

\smallskip
$\bullet$ $s=t$, and $\omega= \sigma_s$; or

\smallskip
$\bullet$ $s=t$, and there exists $r \in S$ such that $m_{s\,r}$ is even and  $\omega= 
w(m_{s\,r}-1: \sigma_s, \sigma_r)$; or

\smallskip
$\bullet$ $s \neq t$, $m_{s\,t}$ is odd, and  $\omega= w(m_{s\,t}-1: \sigma_s, \sigma_t)$.

\smallskip\noindent
Define an {\it elementary $(t,s)$-ribbon} to be either an elementary positive $(t,s)$-ribbon,
or the inverse of an elementary positive $(s,t)$-ribbon.
Note that, if $\omega$ is an elementary $(t,s)$-ribbon, then $\omega \sigma_s = \sigma_t \omega$ 
and $\omega \tau_s = \tau_t \omega$. We say that $\omega$ is a {\it $(t,s)$-ribbon} if there 
exist a sequence $s_0=s, s_1, \dots, s_p=t$ in $S$, and a sequence $\omega_1, \dots, \omega_p$ in 
$A$, such that $\omega_i$ is an elementary $(s_i,s_{i-1})$-ribbon for all $i=1, \dots, p$, and 
$\omega= \omega_p \dots \omega_2 \omega_1$. Clearly, if $\omega$ is a $(t,s)$-ribbon, then 
$\omega \sigma_s = \sigma_t \omega$ and $\omega \tau_s = \tau_t \omega$.

The key point in the proof of Proposition 3.1 is the following result which can be found in 
\cite{God} (see also \cite{God2}).

\begin{prop}[Godelle \cite{God}]
Assume $\Gamma$ to be of type FC. Let $s \in S$, $X \subset S$, and $\alpha \in A$. Then the 
followings are equivalent.

\smallskip
(1) $\alpha \sigma_s \alpha^{-1} \in A_X$.

\smallskip
(2) There exists $k \in \Z \setminus \{0\}$ such that $\alpha \sigma_s^k \alpha^{-1} \in A_X$.

\smallskip
(3) There exist $t \in X$, $\omega \in A$, and $\beta \in A_X$, such that $\omega$ is a $(t,s)$-ribbon 
and $\alpha= \beta \omega$.
\noproof
\end{prop}

\begin{corollary}[Godelle \cite{God}]
Assume $\Gamma$ to be of type FC. Let $s,t \in S$ and $\alpha \in A$. Then the followings are 
equivalent.

\smallskip
(1) $\alpha \sigma_s \alpha^{-1} = \sigma_t$.

\smallskip
(2) There exists $k \in \Z \setminus \{0\}$ such that $\alpha \sigma_s^k \alpha^{-1} = 
\sigma_t^k$.

\smallskip
(3) $\alpha$ is a $(t,s)$-ribbon.
\noproof
\end{corollary}

\begin{lemma}
Assume $\Gamma$ to be of type FC. Then there exists a bijection $\varphi: \hat \Upsilon \to 
\Upsilon$ which sends $\alpha \tau_s \alpha^{-1}$ to $\alpha \sigma_s \alpha^{-1}$ for all 
$\alpha \in A$ and all $s \in S$.
\end{lemma}

\proof
We take $\alpha, \beta \in A$ and $s,t \in S$, and we turn to prove the following equivalence.
$$
\alpha \tau_s \alpha^{-1} = \beta \tau_t \beta^{-1} \quad \Leftrightarrow\quad \alpha \sigma_s 
\alpha^{-1} = \beta \sigma_t \beta^{-1} \,.
$$

Assume first that $\alpha \tau_s \alpha^{-1} = \beta \tau_t \beta^{-1}$. Recall the epimorphism 
$\theta: SA \to A$ which sends $\sigma_s^{\pm 1}$ to $\sigma_s^{\pm 1}$ and $\tau_s$ to 
$\sigma_s$ for all $s \in S$. Then
$$
\alpha \sigma_s \alpha^{-1} = \theta( \alpha \tau_s \alpha^{-1}) = \theta( \beta \tau_t \beta^{-
1}) = \beta \sigma_t \beta^{-1} \,.
$$

Now, assume that $\alpha \sigma_s \alpha^{-1} = \beta \sigma_t \beta^{-1}$. By Corollary 3.4, 
$\beta^{-1} \alpha$ is a $(t,s)$-ribbon, thus $\beta^{-1} \alpha \tau_s= \tau_t \beta^{-1} 
\alpha$, and therefore $\alpha \tau_s \alpha^{-1} = \beta \tau_t \beta^{-1}$.
\endproof

\begin{lemma}
Assume $\Gamma$ to be of spherical type.
Then the bijection $\varphi: \hat \Upsilon \to \Upsilon$ extends to an isomorphism $\varphi: \hat 
\Omega \to \Omega$.
\end{lemma}

\proof
We take $\alpha, \beta \in A$ and $s,t \in S$, and we turn to prove the following equivalence.
$$
(\alpha \tau_s \alpha^{-1}) (\beta \tau_t \beta^{-1}) = (\beta \tau_t \beta^{-1}) (\alpha \tau_s 
\alpha^{-1}) \quad \Leftrightarrow \quad (\alpha \sigma_s \alpha^{-1}) (\beta \sigma_t \beta^{-
1}) = (\beta \sigma_t \beta^{-1}) (\alpha \sigma_s \alpha^{-1})\,.
$$

First, assume that $(\alpha \tau_s \alpha^{-1}) (\beta \tau_t \beta^{-1}) = (\beta \tau_t 
\beta^{-1}) (\alpha \tau_s \alpha^{-1})$. Then
$$
(\alpha \sigma_s \alpha^{-1}) (\beta \sigma_t \beta^{-1}) = \theta( (\alpha \tau_s \alpha^{-1}) 
(\beta \tau_t \beta^{-1})) = \theta( (\beta \tau_t \beta^{-1}) (\alpha \tau_s \alpha^{-1})) = 
(\beta \sigma_t \beta^{-1}) (\alpha \sigma_s \alpha^{-1})\,.
$$

Now, we assume that $(\alpha \sigma_s \alpha^{-1}) (\beta \sigma_t \beta^{-1}) = (\beta \sigma_t 
\beta^{-1}) (\alpha \sigma_s \alpha^{-1})$, and we prove that 
\break
$(\alpha \tau_s \alpha^{-1}) (\beta 
\tau_t \beta^{-1}) = (\beta \tau_t \beta^{-1}) (\alpha \tau_s \alpha^{-1})$. Our proof is divided 
into 3 steps.

\bigskip\noindent
{\it Step 1:} Asssume $\beta=1$, $\alpha=a \in A^+$, and $a \sigma_s \wedge_R a=1$. So, $a \sigma_s a^{-
1}$ is a Charney form. Since $\sigma_t a \sigma_s a^{-1} \sigma_t^{-1} = a \sigma_s a^{-1}$, by 
Proposition 3.2, there exists $c \in A^+$ such that $\sigma_t a = ac$ and $\sigma_t a \sigma_s= a 
\sigma_s c$. The element $c$ is clearly of length 1, namely, $c=\sigma_r$ for some $r \in S$, 
and, by Corollary 3.4, the equality $\sigma_t a = a \sigma_r$ implies that $a$ is a $(t,r)$-ribbon. 
Moreover, we have $a \sigma_r \sigma_s = \sigma_t a \sigma_s = a \sigma_s \sigma_r$, thus 
$\sigma_r \sigma_s = \sigma_s \sigma_r$, therefore $m_{s\,r}=2$. So,
$$
\tau_t a \tau_s a^{-1} = a \tau_r \tau_s a^{-1} = a \tau_s \tau_r a^{-1} = a \tau_s a^{-1} 
\tau_t\,.
$$

\noindent
{\it Step 2:} Assume $\beta=1$ and $\alpha=a \in A^+$. We argue by induction on the length of $a$. The 
case $a \sigma_s \wedge_R a=1$ is treated in Step 1, thus we can suppose that $a \sigma_s 
\wedge_R a \neq 1$. In particular, there exists some $r \in S$ such that $\sigma_r \le_R a 
\sigma_s \wedge_R a$. Suppose $r=s$. Then $a$ can be written as $a=a_1 \sigma_s$ with $a_1 \in 
A^+$, and, moreover, $\sigma_t a_1 \sigma_s a_1^{-1}= a_1 \sigma_s a_1^{-1} \sigma_t$. By the 
inductive hypothesis, it follows that
$$
\tau_t a \tau_s a^{-1} = \tau_t a_1 \tau_s a_1^{-1} = a_1 \tau_s a_1^{-1} \tau_t = a \tau_s a^{-
1} \tau_t\,.
$$
Suppose $r \neq s$. We have $\sigma_s, \sigma_r \le_R a\sigma_s$, thus, by Proposition 3.2, 
$\sigma_s \vee_R \sigma_r = w(m_{s\,r}: \sigma_r, \sigma_s) \le_R a \sigma_s$, therefore $\omega= 
w(m_{s\,r}-1: \sigma_s, \sigma_r) \le_R a$. Write $u=s$ if $m_{s\,r}$ is even, and $u=r$ if 
$m_{s\,r}$ is odd. Then $\omega$ is an (elementary) $(u,s)$-ribbon and $a$ can be written as 
$a=a_1 \omega$ with $a_1 \in A^+$. Moreover, $\sigma_t a_1 \sigma_u a_1^{-1} = \sigma_t a 
\sigma_s a^{-1} = a \sigma_s a^{-1} \sigma_t = a_1 \sigma_u a_1^{-1} \sigma_t$. By the inductive 
hypothesis, it follows that
$$
\tau_t a \tau_s a^{-1} = \tau_t a_1 \tau_u a_1^{-1} = a_1 \tau_u a_1^{-1} \tau_t = a \tau_s a^{-
1} \tau_t\,.
$$

\noindent
{\it Step 3:} General case. By Proposition 3.2, there exist $a \in A^+$ and $k \in \Z$ such that 
$\beta^{-1} \alpha= a \Delta^k$. Recall the permutation $\mu: S \to S$ such that $\Delta \sigma_r 
\Delta^{-1} = \sigma_{\mu(r)}$ for all $r \in S$. Let $r = \mu^k(s)$. Then $\Delta^k$ is a 
$(r,s)$-ribbon (by Corollary 3.4), and, moreover, $\sigma_t a \sigma_r a^{-1} = \sigma_t (\beta^{-1} 
\alpha \sigma_s 
\alpha^{-1} \beta) = (\beta^{-1} \alpha \sigma_s \alpha^{-1} \beta) \sigma_t = a \sigma_r a^{-1} 
\sigma_t$. By Step 2, it follows that
$$
\tau_t (\beta^{-1} \alpha \tau_s \alpha^{-1} \beta) = \tau_t a \tau_r a^{-1} = a \tau_r a^{-1} 
\tau_t = (\beta^{-1} \alpha \tau_s \alpha^{-1} \beta) \tau_s\,,
$$
hence $(\beta \tau_t \beta^{-1})(\alpha \tau_s \alpha^{-1}) = (\alpha \tau_s \alpha^{-1}) (\beta 
\tau_t \beta^{-1})$.
\endproof

\bigskip
In order to extend this result to all Artin groups of type FC (namely, in order to prove 
Proposition 3.1), we need one more preliminary result (on amalgamated products).

\begin{prop}[Serre \cite{Ser}]
Let $G=G_1 \ast_H G_2$ be the amalgamated product of two groups $G_1$ and $G_2$ over $H$. Let 
$C_1$ and $C_2$ be transversals of $G_1/H$ and $G_2/H$, respectively, which contain 1. For all $g 
\in G$ there exists a unique sequence $(g_1, \dots, g_l,h)$ such that

\smallskip
$\bullet$ $g=g_1g_2 \dots g_lh$;

\smallskip
$\bullet$ $h\in H$, and either $g_i \in C_1 \setminus \{1\}$ or $g_i \in C_2 \setminus \{1\}$ for 
all $i=1, \dots, l$;

\smallskip
$\bullet$ $g_i \in C_1 \Rightarrow g_{i+1} \in C_2$, and $g_i \in C_2 \Rightarrow g_{i+1} \in 
C_1$, for all $i=1, \dots, l-1$.
\noproof
\end{prop}

\bigskip
The above sequence $(g_1, \dots, g_l,h)$ is called the {\it amalgam normal form} of $g$ (relative 
to the amalgamed product $G_1 \ast_H G_2$). The number $l$ is called the {\it amalgam norm} of 
$g$ (relative to the amalgamed product $G_1 \ast_H G_2$) and is denoted by $l=|g|_\ast$.
\paragraph{Proof of Proposition 3.1.}
We take $\alpha, \beta \in A$ and $s,t \in S$, and we turn to prove the following equivalence.
$$
(\alpha \tau_s \alpha^{-1}) (\beta \tau_t \beta^{-1}) = (\beta \tau_t \beta^{-1}) (\alpha \tau_s 
\alpha^{-1}) \quad \Leftrightarrow \quad (\alpha \sigma_s \alpha^{-1}) (\beta \sigma_t \beta^{-
1}) = (\beta \sigma_t \beta^{-1}) (\alpha \sigma_s \alpha^{-1})\,.
$$

The implication $\Rightarrow$ can be proved exactly in the same manner as the implication 
$\Rightarrow$ in the proof of Lemma 3.6. So, we assume that $(\alpha \sigma_s \alpha^{-1}) (\beta 
\sigma_t \beta^{-1}) = (\beta \sigma_t \beta^{-1}) (\alpha \sigma_s \alpha^{-1})$, and we prove that 
$(\alpha \tau_s \alpha^{-1}) (\beta \tau_t \beta^{-1}) = (\beta \tau_t \beta^{-1}) (\alpha \tau_s 
\alpha^{-1})$. We argue by induction on the rank of $A$ ({\it i.e.} on the number of vertices of 
$\Gamma$). Up to changing $\alpha$ by $\beta^{-1} \alpha$, we can also assume that $\beta=1$.

If $m_{x\,y} <+\infty$ for all $x,y \in S$, $x\neq y$, then $A$ is of spherical type and, 
therefore, the equality $(\alpha \tau_s \alpha^{-1}) \tau_t = \tau_t (\alpha \tau_s \alpha^{-1})$ 
follows from Lemma 3.6. So, we can assume that there exist $x,y \in S$, $x \neq y$, such that 
$m_{x\,y}= +\infty$. Let $X=S \setminus \{x\}$ and $Y=S \setminus \{y\}$. Then $A=A_X \ast_{A_{X 
\cap Y}} A_Y$. 
We choose transversals $C_X$ and $C_Y$ of $A_X/A_{X\cap Y}$ and $A_Y/A_{X \cap Y}$, respectively,
which contain 1,
we consider the amalgam normal form $(\alpha_1, \dots, \alpha_l, \gamma)$ of 
$\alpha$, and we argue by induction on $l=|\alpha|_\ast$.

Assume $l=0$. So, $\alpha=\gamma \in A_{X \cap Y}$. The fact that the amalgam normal forms 
of $\sigma_t (\alpha \sigma_s \alpha^{-1})$ and $(\alpha \sigma_s \alpha^{-1}) \sigma_t$ are 
equal implies that either $\sigma_s, \sigma_t \in A_X$ (namely $s,t \in X$), or $\sigma_s, 
\sigma_t \in A_Y$ (namely, $s,t \in Y$). Say $s,t \in X$. We have $s,t \in X$ and $\alpha \in 
A_X$, thus, by the inductive hypothesis, $\tau_t (\alpha \tau_s \alpha^{-1}) = (\alpha \tau_s 
\alpha^{-1}) \tau_t$.

Now, we assume that $l=|\alpha|_\ast \ge 1$. Without loss of generality, we can assume that 
$\alpha_l \in C_X$. Then $\alpha_1 \in C_Z$, where $Z=Y$ if $l$ is even, and $Z=X$ if $l$ is odd. 
We consider 4 different cases.

\bigskip\noindent
{\it Case 1:} $\alpha_l \gamma \sigma_s \gamma^{-1} \alpha_l^{-1} \in A_{X \cap Y}$. By 
Proposition 3.3, $\alpha_l \gamma$ can be written as $\alpha_l \gamma = \gamma_1 \omega$ where 
$\omega$ is a $(r,s)$-ribbon for some $r \in X \cap Y$, and $\gamma_1 \in A_{X \cap Y}$. Let 
$\alpha'= \alpha_1 \alpha_2 \dots \alpha_{l-1} \gamma_1= \alpha \omega^{-1}$. We have $(\alpha' 
\sigma_r {\alpha'}^{-1}) \sigma_t = (\alpha \sigma_s \alpha^{-1}) \sigma_t = \sigma_t (\alpha 
\sigma_s \alpha^{-1}) = \sigma_t (\alpha' \sigma_r {\alpha'}^{-1})$, thus, by the inductive 
hypothesis (on the amalgam norm of $\alpha$), we have $(\alpha \tau_s \alpha^{-1}) \tau_t = 
(\alpha' \tau_r {\alpha'}^{-1}) \tau_t = \tau_t (\alpha' \tau_r {\alpha'}^{-1}) = \tau_t (\alpha 
\tau_s \alpha^{-1})$.

\bigskip\noindent
{\it Case 2:} $l=1$, $s \in X$, and $\alpha \sigma_s \alpha^{-1} \in A_X \setminus A_{X\cap Y}$. 
The fact that the amalgam normal forms of $\sigma_t (\alpha \sigma_s \alpha^{-1})$ 
and $(\alpha \sigma_s \alpha^{-1}) \sigma_t$ are equal implies that $\sigma_t \in A_X$, namely, that $t 
\in X$. We have $s,t \in X$ and $\alpha \in A_X$, thus, by the inductive hypothesis (on the rank 
of $A$), we have $\tau_t (\alpha \tau_s \alpha^{-1})= (\alpha \tau_s \alpha^{-1}) \tau_t$.

\bigskip\noindent
{\it Case 3:} $l\ge 2$, $s \in X$, and $\alpha_l \gamma \sigma_s \gamma^{-1} \alpha_l^{-1} \in 
A_X \setminus A_{X \cap Y}$. Then the amalgam normal form of $\alpha \sigma_s \alpha^{-1}$ has 
the form $(\alpha_1, \dots, \alpha_{l-1}, \beta_1', \beta_2', \dots, \beta_l', \gamma_1)$, and 
$\alpha_1,\beta_l' \in C_Z$. The fact that the amalgam normal forms of 
$\sigma_t (\alpha \sigma_s \alpha^{-1})$ and $(\alpha \sigma_s \alpha^{-1}) \sigma_t$ are equal 
implies that $\sigma_t \in A_Z$, namely, that $t \in Z$. Then the amalgam normal form of $(\alpha \sigma_s 
\alpha^{-1}) \sigma_t$ has either the form $(\alpha_1, \dots, \alpha_{l-1}, \beta_1', \beta_2', 
\dots, \beta_{l-1}', \gamma_2)$ if $\beta_l' \gamma_1 \sigma_t \in A_{X \cap Y}$, or the form 
$(\alpha_1, \dots, \alpha_{l-1}, 
\break
\beta_1', \beta_2', \dots, \beta_{l-1}', \beta_l'', \gamma_2)$ if 
$\beta_l' \gamma_1 \sigma_t \not\in A_{X \cap Y}$. Now, $\sigma_t (\alpha \sigma_s \alpha^{-1})$ 
has also this amalgam normal form, thus $\alpha_1 A_{X \cap Y} = \sigma_t \alpha_1 A_{X\cap Y}$, 
namely, $\alpha_1^{-1} \sigma_t \alpha_1 \in A_{X\cap Y}$. By Proposition 3.3, we deduce that 
$\alpha_1$ can be written as $\alpha_1=\omega \delta$, where $\omega$ is a $(t,r)$-ribbon for 
some $r \in X \cap Y$, and $\delta \in A_{X \cap Y}$. Let $\alpha'= \delta \alpha_2 \dots 
\alpha_l \gamma = \omega^{-1} \alpha$. We have $\sigma_r (\alpha' \sigma_s {\alpha'}^{-1}) = 
(\alpha' \sigma_s {\alpha'}^{-1}) \sigma_r$, thus, by the inductive hypothesis (on the amalgam 
norm of $\alpha$), we have $\tau_r (\alpha' \tau_s {\alpha'}^{-1}) = (\alpha' \tau_s {\alpha'}^{-
1}) \tau_r$, therefore $\tau_t (\alpha \tau_s \alpha^{-1}) = (\alpha \tau_s \alpha^{-1}) \tau_t$.

\bigskip\noindent
{\it Case 4:} $l \ge 1$ and $s=x \not\in X$. Then the amalgam normal form of $\alpha \sigma_s 
\alpha^{-1}$ has the form $(\alpha_1, \dots, \alpha_{l-1}, \alpha_l, \beta_0', \beta_1', \dots, 
\beta_l', \gamma_1)$. Applying the same argument as in Case 3, we conclude that $\tau_t (\alpha 
\tau_s \alpha^{-1}) = (\alpha \tau_s \alpha^{-1}) \tau_t$.
\endproof

\section{The FRZ property}

The aim of this section is to prove the following.

\begin{prop}
Assume $\Gamma$ to be of type FC. Let $\alpha \in SA$ and $s,t \in S$. Then the followings are 
equivalent.

\smallskip
(1) $\alpha \sigma_s = \sigma_t \alpha$.

\smallskip
(2) There exists $k \in \Z \setminus \{0\}$ such that $\alpha \sigma_s^k = \sigma_t^k \alpha$.

\smallskip
(3) $\alpha \tau_s = \tau_t \alpha$.

\smallskip
(4) There exists $k \in \N \setminus \{0\}$ such that $\alpha \tau_s^k = \tau_t^k \alpha$.
\end{prop}

The following lemma is a preliminary result to the proof of Proposition 4.1.

\begin{lemma}
Let $G$ be a graph, let $V$ be its set of vertices, and let $E$ be its set of edges. Let $u_1, \dots, u_l,v,w 
\in V$ and $p \in \N \setminus \{0\}$ such that $v^p u_1 u_2 \dots u_l = u_1 u_2 \dots u_l w^p$ 
(in $\MM (G)$). Then $v=w$, and $u_iv=vu_i$ for all $i=1, \dots, l$.
\end{lemma}

\proof
Observe that, if $H$ is a full subgraph of $G$, then there is an epimorphism $f_H: \MM(H) \to
\MM(G)$ which sends the vertices of $H$ to themselves and sends the other vertices of $G$ to $1$.  
First, applying this observation to the graph $H=\{v\}$ with one vertex, $v$, and
no edge, we deduce that $v=w$. Now, we argue by induction on $l$. Suppose that $u_1 \neq v$ and that
$u_1$ and $v$ are not joined by an edge. Let $H$ be the full subgraph of $G$ generated by $\{u_1,v\}$.
Then $\MM (H)$ is the free monoid freely generated by $\{u_1,v\}$, $f_M(v^pu_1 \dots u_l)$ is a
word which starts with $v$, and the word $f_M(u_1 \dots u_l v^p)$ starts with $u_1$: a contradiction.
So, either $u_1=v$ or $\{u_1,v\} \in E$, that is, $u_1 v=v u_1$. It follows also that
$v^p u_2 \dots u_l = u_2 \dots u_l v^p$, and, by the inductive hypothesis, we conclude that
$u_i v=v u_i$ for all $i=2, \dots, l$.
\endproof

\paragraph{Proof of Proposition 4.1.}
The implications (1)$\Rightarrow$(2) and (3)$\Rightarrow$(4) are obvious, thus it remains to 
prove (2)$\Rightarrow$(3) and (4)$\Rightarrow$(1).

\bigskip\noindent
{\it Proof of (2)$\Rightarrow$(3).}
We assume that $\alpha \sigma_s^k = \sigma_t^k \alpha$ for some $k \in \Z \setminus \{0\}$. 
Recall the epimorphism $\theta: SA \to A$ which sends $\sigma_s^{\pm 1}$ to $\sigma_s^{\pm 1}$ 
and $\tau_s$ to $\sigma_s$ for all $s \in S$. Let $\alpha_0= \theta (\alpha) \in A$. Then 
$\alpha_0 \sigma_s^k = \sigma_t^k \alpha_0$, thus, by Corollary 3.4, $\alpha_0$ is a $(t,s)$-ribbon, 
therefore $\alpha_0 \tau_s = \tau_t 
\alpha_0$. Let $\beta= \alpha_0^{-1} \alpha$. Then $\beta \sigma_s^k = \sigma_s^k \beta$. Recall 
the decomposition $SA= \MM( \hat \Omega) \rtimes A$, and write $\beta = u_1 u_2 \dots u_l 
\beta_0$ where $u_1, \dots, u_l \in \hat \Upsilon$ and $\beta_0 \in A$. We have $\sigma_s^k \beta 
\sigma_s^{-k}= u_1' u_2' \dots u_l' \beta_0'$, where $u_i'= \sigma_s^k u_i \sigma_s^{-k} \in \hat 
\Upsilon$ for all $i=1, \dots, l$, and $\beta_0'= \sigma_s^k \beta_0 \sigma_s^{-k} \in A$. The 
equality $\sigma_s^k \beta \sigma_s^{-k} = \beta$ implies that $\beta_0'=\beta_0$, and that there 
exists a permutation $\chi\in \Sym_l$ such that $u_i'= u_{\chi(i)}$ for all $i=1, \dots, l$. 
Firstly, the equality $\beta_0'= \sigma_s^k \beta_0 \sigma_s^{-k} = \beta_0$ implies by Corollary 
3.4 that $\beta_0$ is a $(s,s)$-ribbon, thus $\beta_0 \tau_s = \tau_s \beta_0$. Now, let $p$ be 
the order of $\chi$. Then $\sigma_s^{pk} u_i \sigma_s^{-pk} = u_i$ for all $i=1, \dots, l$. Take 
$i\in \{1, \dots, l\}$ and write $u_i=\gamma_i \tau_{r_i} \gamma_i^{-1}$ where $\gamma_i \in A$ and 
$r_i \in S$. The element $\sigma_s^{pk}$ commutes with $\gamma_i \tau_{r_i} \gamma_i^{-1} = u_i$, 
thus $\sigma_s^{pk}$ commutes with $\gamma_i \sigma_{r_i} \gamma_i^{-1} = \theta( \gamma_i 
\tau_{r_i} \gamma_i^{-1})$, therefore $\sigma_s$ commutes with $\gamma_i \sigma_{r_i} \gamma_i^{-
1}$ (by Corollary 3.4), and hence, by Proposition 3.1, $\tau_s$ commutes with $\gamma_i 
\tau_{r_i} \gamma_i^{-1} = u_i$. This shows that $\tau_s$ commutes with $\beta= u_1 \dots u_l 
\beta_0$, and we conclude that $\alpha \tau_s = \alpha_0 \beta \tau_s = \alpha_0 \tau_s \beta = 
\tau_t \alpha_0 \beta = \tau_t \alpha$.

\bigskip\noindent
{\it Proof of (4)$\Rightarrow$(1).} 
We assume that $\alpha \tau_s^k = \tau_t^k \alpha$ for some $k \in \N \setminus \{0\}$. We write 
$\alpha= u_1 u_2 \dots u_l \alpha_0$, where $u_1, \dots, u_l \in \hat \Upsilon$ and $\alpha_0 \in 
A$. Then the equality $\alpha \tau_s^k = \tau_t^k \alpha$ implies that we have $u_1 
u_2 \dots u_l (\alpha_0 \tau_s \alpha_0^{-1})^k = \tau_t^k u_1 u_2 \dots u_l$ in $\MM (\hat \Omega)$. 
By Lemma~4.2, it follows that $\alpha_0 \tau_s \alpha_0^{-1} = \tau_t$, and $\tau_t u_i = u_i 
\tau_t$ for all $i=1, \dots, l$. The equality $\alpha_0 \tau_s \alpha_0^{-1} = \tau_t$ implies
that $\alpha_0 \sigma_s \alpha_0^{-1} = \theta (\alpha_0 \tau_s \alpha_0^{-1}) = \theta(\tau_t)
= \sigma_t$. Now, take $i\in \{1, \dots, 
l\}$ and write $u_i= \gamma_i \tau_{r_i} \gamma_i^{-1}$ where $\gamma_i \in A$ and $r_i \in S$. 
The element $\tau_t$ commutes with $u_i=\gamma_i \tau_{r_i} \gamma_i^{-1}$, thus $\sigma_t= 
\theta (\tau_t)$ commutes with $\gamma_i \sigma_{r_i} \gamma_i^{-1}= \theta(u_i)$, that is, 
$\sigma_t \gamma_i \sigma_{r_i} \gamma_i^{-1} \sigma_t^{-1} = \gamma_i \sigma_{r_i} \gamma_i^{-
1}$, hence, by Lemma~3.5, $\sigma_t u_i \sigma_t^{-1} = \sigma_t \gamma_i \tau_{r_i} 
\gamma_i^{-1} \sigma_t^{-1} = \gamma_i \tau_{r_i} \gamma_i^{-1} = u_i$. This shows that $\alpha 
\sigma_s = u_1 \dots u_l \alpha_0 \sigma_s = u_1 \dots u_l \sigma_t \alpha_0 = \sigma_t u_1 \dots 
u_l \alpha_0 = \sigma_t \alpha$.
\endproof 

\section{Birman's conjecture}

The purpose of this section is to prove Birman's conjecture 
for right-angled Artin groups. As these groups are the same as graph groups, for convenience, we 
shall use the terminology of graph groups and graph products of groups.

Let $G$ be a graph, let $V$ be its set of vertices, and let $E$ be its set of edges. Take an 
(abstract) set $\SS= \{ \sigma_u; u \in V\}$ in one-to-one correspondence with $V$, and define 
the {\it graph group} associated to $G$ to be the group presented by
$$
\GG(G)= \langle \SS\ |\ \sigma_u \sigma_v = \sigma_v \sigma_u \text{ for } \{u,v\} \in E 
\rangle\,.
$$
Take another abstract set $\TT= \{ \tau_u; u\in V\}$ in one-to-one correspondence with $V$, and 
define the {\it singular graph monoid} associated to $G$ to be the monoid $S\GG (G)$ 
presented as a monoid by the generating set
$\SS \cup \SS^{-1} \cup \TT$ and by the relations
\[
\begin{array}{cl}
\sigma_u \sigma_u^{-1} = \sigma_u^{-1} \sigma_u = 1\,, \quad \sigma_u \tau_u = \tau_u \sigma_u\,, 
\quad &\text{for } u \in V\,,\\
\sigma_u \sigma_v = \sigma_v \sigma_u\,, \quad \sigma_u \tau_v = \tau_v \sigma_u\,, \quad \tau_u 
\tau_v = \tau_v \tau_u\,, \quad &\text{for } \{u,v\} \in E\,.\\
\end{array}
\]
The {\it desingularization map} is defined in this context as the multiplicative homomorphism 
$\eta: S\GG (G) \to \Z[ \GG (G)]$ which sends $\sigma_u^{\pm 1}$ to $\sigma_u^{\pm 1}$ and 
$\tau_u$ to $\sigma_u - \sigma_u^{-1}$ for all $u \in V$.

\begin{thm}
The desingularization map $\eta: S\GG (G) \to \Z [\GG (G)]$ is injective.
\end{thm}

We begin with some definitions and results on graph products of groups.

Suppose given a group (or a monoid) $K_u$ for each $u \in V$. Then the {\it graph product} of the 
family $\{K_u\}_{u \in V}$ along the graph $G$ is the quotient $\GG( \{K_u\}_{u \in V}, G)$ of the 
free product $\ast_{u \in V} K_u$ by the relations
$$
g_u g_v = g_v g_u\quad \text{for } \{u,v\} \in E,\ g_u \in K_u,\ g_v \in K_v\,.
$$
Note that, if $K_u=\N$ for all $u \in V$, then $\GG(\{K_u\}, G)=\MM (G)$, if $K_u=\Z$ for all $u 
\in V$, then $\GG(\{K_u\},G)= \GG (G)$, and if $K_u= \Z \times \N$ for all $u \in V$, then $\GG 
(\{K_u\},G)= S\GG (G)$.

Let $\GG (\{K_u\},G)$ be a graph product of groups (or monoids). Let $g \in \GG (\{K_u\},G)$. 
An {\it expression} for $g$ is a sequence $W=(g_1, g_2, \dots, g_l)$ such that

\smallskip
$\bullet$ there exists $u_i \in V$ such that $g_i \in K_{u_i} \setminus \{1\}$ for all $i=1, 
\dots, l$,

\smallskip
$\bullet$ $g=g_1g_2 \dots g_l$.

\smallskip\noindent
The {\it support} of $W$ is the sequence $\Supp(W)= (u_1, u_2, \dots, u_l)$, and the {\it length} 
of $W$ is $|W|=l$. The minimal length for an expression for $g$ is called the {\it syllable 
length} of $g$ and is denoted by $|g|_G$.
An expression $W$ of $g$ is called {\it reduced} if its length is equal to the length of $g$.

Let $g \in \GG (\{K_u\},G)$, let $W=(g_1, \dots, g_l)$ be an expression for $g$, and let $(u_1, 
u_2, \dots, u_l)$ be the support of $W$. Suppose there exists $i \in \{1, \dots, l-1\}$ such 
that $u_i=u_{i+1}$, and put
$$
W'=\left\{
\begin{array}{ll}
(\dots, g_{i-1}, g_{i+2}, \dots) \quad &\text{if } g_ig_{i+1}=1\,,\\
(\dots, g_{i-1}, g_ig_{i+1}, g_{i+2}, \dots) \quad &\text{if } g_ig_{i+1} \neq 1\,.\\
\end{array}
\right.
$$
We say that $W'$ is obtained from $W$ via an {\it elementary M-operation of type I}. This 
operation shortens the length of an expression by 1 or 2. Suppose there exists $i \in \{ 1, \dots, 
l-1\}$ such that $\{u_i,u_{i+1} \} \in E$, and put
$$
W''= (\dots, g_{i-1}, g_{i+1}, g_i, g_{i+2}, \dots)\,.
$$
We say that $W''$ is obtained from $W$ via an {\it elementary M-operation of type II}. This 
operation leaves the length of an expression unchanged. We say that an expression $W$ is {\it M-
reduced} if its length cannot be reduced applying a sequence of elementary M-operations.

The following proposition is essentially proved in \cite{Gre} (see also \cite{HeMe} and 
\cite{HsWi}).

\begin{prop}[Green \cite{Gre}]
(1) Let $g \in \GG(\{K_u\},G)$, and let $W_1$ and $W_2$ be two M-reduced expressions for $g$. 
Then $W_1$ and $W_2$ are related by a sequence of elementary M-operations of type II.

\smallskip
(2) Let $g \in \GG(\{K_u\},G)$, and let $W$ be an expression for $g$. Then $W$ is M-reduced if 
and only if $W$ is reduced.
\noproof
\end{prop}

Now, assume that the set $V$ of vertices is endowed with a total order, $\le$. Let $g \in 
\GG(\{K_u\},G)$, and let $W$ be a reduced expression for $g$. We say that $W$ is a {\it normal 
form} for $g$ if $\Supp (W)$ is the smallest support of a reduced expression for $g$ with respect 
to the lexicographic order.

\begin{corollary}
(1) Each element of $\GG(\{K_u\},G)$ has a unique normal form.

\smallskip
(2) To be a normal form depends only on its support, namely, if $W$ is the normal form for some
$g \in \GG(\{K_u\},G)$, if $W'$ is some expression for some $g' \in \GG(\{K_u\},G)$, and if
$\Supp (W')= \Supp (W)$, then $W'$ is the normal form for $g'$.
\noproof
\end{corollary}

We return to the study of the graph group $\GG (G)$. Consider the homomorphism $\deg: \GG (G) 
\to \Z$ which sends $\sigma_u$ to $1$ for all $u \in V$. For $n \in \Z$, we put $\GG_n (G)= \{ 
g\in \GG(G); \deg(g) \ge n\}$, and we denote by $\AA_n$ the free $\Z$-module freely generated by 
$\GG_n (G)$. Let $\AA= \Z [\GG (G)]$. Then $\{ \AA_n\}_{n \in \Z}$ is a filtration of $\AA$ 
compatible with the multiplication, that is,

\smallskip
$\bullet$ $\AA_n \subset \AA_m$ if $n \ge m$,

\smallskip
$\bullet$ $\AA_p \AA_q \subset \AA_{p+q}$ for $p,q \in \Z$,

\smallskip
$\bullet$ $1 \in \AA_0$.

\smallskip\noindent
Moreover, this filtration is separate, namely,

\smallskip
$\bullet$ $\cap_{n \in \Z} \AA_n = \{0\}$.

\smallskip\noindent
We denote by $\tilde \AA$ the completion of $\AA$. For $n \in \Z$, we put $\GG^{(n)}(G)= \{g \in 
\GG (G); \deg(g)=n\}$, and we denote by $\AA^{(n)}$ the free $\Z$-module freely generated by 
$\GG^{(n)} (G)$. Then each element of $\tilde \AA$ can be uniquely written as a formal series 
$\sum_{n=k}^{+\infty} P_n$, where $k \in \Z$ (which may be negative), and $P_n \in \AA^{(n)}$ for 
all $n \ge k$.

Take a free abelian group $K_u \simeq \Z \times \Z$ of rank $2$ generated by $\{\sigma_u, 
\tau_u\}$ for all $u \in V$, and write $\tilde \GG= \GG(\{K_u\},G)$. Let $\UU (\tilde \AA)$ be the 
group of unities of $\tilde \AA$. We have a homomorphism $\tilde \eta: \tilde \GG \to \UU( \tilde 
\AA)$ defined by
$$
\tilde \eta (\sigma_u)= \sigma_u\,, \quad \tilde \eta (\tau_u)= \sigma_u - \sigma_u^{-1}\,, \quad 
\text{for } u\in V\,.
$$
Note that
$$
\tilde \eta (\tau_u^{-1}) = -\sum_{n=0}^{+\infty} \sigma_u^{2n+1}\,, \quad \text{for } u \in V\,.
$$
Note also that $S\GG (G)$ is a submonoid of $\tilde \GG$, that $\AA= \Z [\GG (G)]$ is a 
subalgebra of $\tilde \AA$, and that 
the restriction of $\tilde \eta$ to $S\GG (G)$ is the desingularization map $\eta: S\GG(G) \to \Z 
[\GG (G)]$. So, Theorem 5.1 is a consequence of the following.  

\begin{prop}
The homomorphism $\tilde \eta: \tilde \GG \to \UU( \tilde \AA)$ is injective.
\end{prop}

\proof
For $u \in V$ and $(p,q) \in \Z \times \Z$, $(p,q) \neq (0,0)$, we write
$$
(\sigma_u - \sigma_u^{-1})^p \sigma_u^q = \sum_{n=k}^{+\infty} c_{n\,p\,q} \sigma_u^n\,.
$$
Note that the numbers $c_{n\,p\,q}$ do not depend on $u$, but only on $n$, $p$, 
and $q$. Note also that there always exists some $a \ge k$ (which may be negative) such that $a \neq 0$ 
and $C_{a\,p\,q} 
\neq 0$.

Let $\tilde g \in \tilde \GG$, $\tilde g \neq 1$. Let $(\tau_{u_1}^{p_1} \sigma_{u_1}^{q_1}, 
\tau_{u_2}^{p_2} \sigma_{u_2}^{q_2}, \dots, \tau_{u_l}^{p_l} \sigma_{u_l}^{q_l})$ be the normal 
form for $\tilde g$. We have
$$
\begin{array}{rl}
\tilde \eta(\tilde g) = &(\sigma_{u_1} - \sigma_{u_1}^{-1})^{p_1} \sigma_{u_1}^{q_1} 
(\sigma_{u_2} - \sigma_{u_2}^{-1})^{p_2} \sigma_{u_2}^{q_2} \dots (\sigma_{u_l} - \sigma_{u_l}^{-
1})^{p_l} \sigma_{u_l}^{q_l}\\
= &{\displaystyle \sum_{n_1 \ge k_1, \dots, n_l \ge k_l}} c_{n_1\,p_1\,q_1} c_{n_2\,p_2\,q_2} 
\dots c_{n_l\,p_l\,q_l}\, \sigma_{u_1}^{n_1} \sigma_{u_2}^{n_2} \dots \sigma_{u_l}^{n_l}\,.
\end{array}
$$ 
By the above observations, we can find $a_1, \dots, a_l \in \Z \setminus \{0\}$ such that 
$c_{a_i\,p_i\,q_i} \neq 0$ for all $i=1, \dots, l$. We turn to show that $\sigma_{u_1}^{n_1} 
\sigma_{u_2}^{n_2} \dots \sigma_{u_l}^{n_l} \neq \sigma_{u_1}^{a_1} \sigma_{u_2}^{a_2} \dots 
\sigma_{u_l}^{a_l}$ for any $l$-tulpe $(n_1, \dots, n_l)$ in $\Z^l$ different from $(a_1, \dots, a_l)$. 
This implies that 
$\sigma_{u_1}^{a_1} \sigma_{u_2}^{a_2} \dots \sigma_{u_l}^{a_l} \neq 1$, and that the coefficient 
of $\sigma_{u_1}^{a_1} \sigma_{u_2}^{a_2} \dots \sigma_{u_l}^{a_l}$ in $\tilde \eta 
(\tilde g)$ is $c_{a_1\,p_1\,q_1} c_{a_2\,p_2\,q_2} \dots c_{a_l\,p_l\,q_l} \neq 0$, thus that 
$\tilde \eta (\tilde g) \neq 1$.

The sequence $(\tau_{u_1}^{p_1} \sigma_{u_1}^{q_1}, \dots, \tau_{u_l}^{p_l} \sigma_{u_l}^{q_l})$ 
is a normal form, thus, by Corollary 5.3, the sequence $(\sigma_{u_1}^{a_1}, \sigma_{u_2}^{a_2}, 
\dots, \sigma_{u_l}^{a_l})$ is the normal form for $\sigma_{u_1}^{a_1} \sigma_{u_2}^{a_2} \dots 
\sigma_{u_l}^{a_l}$. Assume $n_i \neq 0$ for all $i=1, \dots, l$. Again, by Corollary 5.3, the 
sequence $(\sigma_{u_1}^{n_1}, \sigma_{u_2}^{n_2}, \dots, \sigma_{u_l}^{n_l})$ is the normal form 
for $\sigma_{u_1}^{n_1} \sigma_{u_2}^{n_2} \dots \sigma_{u_l}^{n_l}$, thus 
$\sigma_{u_1}^{n_1} \sigma_{u_2}^{n_2} \dots \sigma_{u_l}^{n_l} \neq \sigma_{u_1}^{a_1} 
\sigma_{u_2}^{a_2} \dots \sigma_{u_l}^{a_l}$ if $(n_1,n_2, \dots, n_l) \neq (a_1, a_2, \dots, 
a_l)$. Now, assume that there exists $i \in \{1, \dots, l\}$ such that $n_i=0$. Then
$$
|\sigma_{u_1}^{n_1} \sigma_{u_2}^{n_2} \dots \sigma_{u_l}^{n_l} |_G < l =|( \sigma_{u_1}^{a_1}, 
\sigma_{u_2}^{a_2}, \dots, \sigma_{u_l}^{a_l})| = |\sigma_{u_1}^{a_1} \sigma_{u_2}^{a_2} \dots 
\sigma_{u_l}^{a_l}|_G\,,
$$
hence $\sigma_{u_1}^{n_1} \sigma_{u_2}^{n_2} \dots \sigma_{u_l}^{n_l} \neq \sigma_{u_1}^{a_1} 
\sigma_{u_2}^{a_2} \dots \sigma_{u_l}^{a_l}$.
\endproof

\section{Other questions}

\subsection{Topological interpretation}

A {\it geometric braid} on $n$ strings is a $n$-tuple $\beta= (b_1, \dots, b_n)$ of smooth 
disjoint paths in $\C \times [0,1]$ such that

\smallskip
$\bullet$ the projection of $b_k(t)$ on the second component is $t$ for all $k \in \{1, \dots, 
n\}$ and all $t \in [0,1]$,

\smallskip
$\bullet$ $b_k(0)=(k,0)$ and $b_k(1)=(\xi(k),1)$, where $\xi$ is some permutation of $\{1, \dots, 
n\}$, for all $k \in \{1, \dots, n\}$.

\smallskip\noindent
The {\it braid group} on $n$ strings, denoted by $\BB_n$, is the group of isotopy classes of 
braids. Define a {\it singular braid} on $n$ strings as a $n$-tuple $\beta= (b_1, \dots, b_n)$ 
of smooth paths in $\C \times [0,1]$ such that

\smallskip
$\bullet$ the projection of $b_k(t)$ on the second component is $t$ for all $k \in \{1, \dots, 
n\}$ and all $t \in [0,1]$,

\smallskip
$\bullet$ $b_k(0)=(k,0)$ and $b_k(1)=(\xi(k),1)$, where $\xi$ is some permutation of $\{1, \dots, 
n\}$, for all $k \in \{1, \dots, n\}$,

\smallskip
$\bullet$ the strings intersect transversely in finitely many double points.

\smallskip\noindent
The {\it singular braid monoid}, denoted by $S\BB_n$, is the monoid of isotopy classes of 
singular braids.
By \cite{Art} (see also \cite{Art2}), the braid group $\BB_n$ is isomorphic to the Artin group 
$A=A_\Gamma$ associated to the Coxeter graph $\Gamma= A_{n-1}$, and, by \cite{Bir}, $S\BB_n$ is 
the singular Artin monoid associated to $\Gamma=A_{n-1}$. 

Some other Artin groups can be viewed 
as ``geometric braid groups". For instance, the Artin group associated to the Coxeter graph 
$\Gamma=B_n$ is the braid group on $n$ strings of the annulus, and the Artin group associated to 
the graph $\Gamma=D_n$ is an index 2 subgroup of the braid group on $n$ strings of the plane 
endowed with a singular point of degree 2 (see \cite{All}). One can easily verify using the 
techniques of \cite{Gon} that the singular braid monoid of the annulus coincides with the 
singular Artin monoid associated to $\Gamma=B_n$. However, we did not success to find any 
embedding of the singular Artin monoid associated to $\Gamma=D_n$ into the singular braid monoid 
of the plane endowed with a singular point of degree 2, and we suspect that such an embedding 
does not exist.

In other respects, the spherical type Artin groups can be interpreted as fundamental groups of 
regular orbit spaces (see \cite{Bri}), and the groups of type $A_n$ and $D_n$ as geometric 
monodromy groups of simple singularities (see \cite{PeVa}). In both cases, we do not know whether these 
topological interpretations can be extended to the singular Artin monoids, even for the singular 
Artin monoids of type $A_n$ (namely, the singular braid monoids).

\subsection{Vassilev invariants}

An {\it invariant} on $A$ is a set-map $v: A \to H$, where $H$ is some abelian group, or, 
equivalently, a homomorphism $v: \Z [A] \to H$ of $\Z$-modules. Consider the homomorphism $\ord: 
SA \to \N$ which sends $\sigma_s^{\pm 1}$ to $0$ and $\tau_s$ to $1$ for all $s \in S$. For $d 
\in \N$, we put $S_dA= \{ \alpha \in SA; \ord (\alpha)=d \}$. So, $S_dA$ is the set of elements 
of $SA$ that have exactly $d$ ``singularities''. Recall the desingularization map $\eta: SA 
\to \Z [A]$ which sends $\sigma_s^{\pm 1}$ to $\sigma_s^{\pm 1}$ and $\tau_s$ to $\sigma_s - 
\sigma_s^{-1}$ for all $s \in S$. Then define a {\it Vassiliev invariant} of type $d$ as an 
invariant $v: \Z[A] \to H$ which vanishes on $\eta (S_{d+1}A)$.

The following questions have been solved for braid groups (see \cite{Bar}, \cite{Koh}, and 
\cite{Pap}), and remain open for the other Artin groups.

\begin{question}
Do Vassiliev invariants separate the elements of $A$? In other words, given two elements $\alpha, 
\beta \in A$, $\alpha \neq \beta$, does there exist a Vassiliev invariant $v: \Z [A] \to H$ such 
that $v(\alpha) \neq v(\beta)$? 
\end{question}

By a {\it universal Vassiliev invariant} we mean an invariant $Z: \Z [A] \to \AA$ such that, if 
$v: \Z [A] \to H$ is some Vassiliev invariant, then there exists a homomorphism $u: \AA \to H$ of 
$\Z$-modules such that $v=u \circ Z$.

\begin{question}
Describe a universal Vassiliev invariant for $A$ in terms of generators and relations.
\end{question}

\subsection{Conjugacy problem}

Let $\MM$ be a monoid, and let $\GG$ be the group of unities of $\MM$. We say that two elements 
$\alpha, \beta \in \MM$ are {\it conjugate} if there exists some $\gamma \in \GG$ such that 
$\gamma \alpha \gamma^{-1} = \beta$. A {\it solution to the conjugacy problem} in $\MM$ is an 
algorithm which decides whether two given elements $\alpha, \beta \in \MM$ are conjugate or not.

As for the word problem, a solution to the conjugacy problem for $SA$ will give a solution to the 
conjugacy problem for $A$ (which, by the way, is the group of unities of $SA$). So, a reasonable 
approach is to study the conjugacy problem for those singular Artin monoids whose associated 
Artin groups have known solutions to the conjugacy problem. For instance, a solution to the 
conjugacy problem for Artin groups of spherical type can be found in \cite{BrSa} (see also 
\cite{Pic} and \cite{FrGo}).

A solution to the conjugacy problem for singular braid monoids is given in \cite{Ver}. We 
suspect that this algorithm can be extended to all spherical type singular Artin monoids. On the 
other hand, we do not know whether the techniques introduced in the present paper, more 
specifically, Propositions 2.1 and 3.1, can be also used to solve the conjugacy problem in $SA$.

On other respects,
let $\tilde A= \GG (\hat \Omega) \rtimes A$ be the smallest group which contains $SA= \MM 
(\hat \Omega) \rtimes A$. As far as we know, no solution to the conjugacy problem for $\tilde A$ 
is known, even in the case where $SA$ is the singular braid monoid.


\bigskip\bigskip\noindent
\halign{#\hfill\hskip2truecm&#\hfill\cr
Eddy Godelle&Luis Paris\cr
Laboratoire de Math\'ematiques
&Institut de Math\'ematiques\cr
Nicolas Oresme&de Bourgogne\cr
Universit\'e de Caen
&Universit\'e de Bourgogne\cr
UMR 6139 du CNRS, BP 5186
&UMR 5584 du CNRS, BP 47870\cr
14032 Caen cedex
&21078 Dijon cedex\cr
France&France\cr
\noalign{\smallskip}
\texttt{eddy.godelle@math.unicaen.fr}
&\texttt{lparis@u-bourgogne.fr}\cr}

\end{document}